\documentclass[11pt,reqno]{amsart}

\usepackage{graphicx}
\usepackage{amsmath,amsthm,amscd,amssymb}
\usepackage{amsfonts}
\usepackage{latexsym}

\newcommand{\R}{\mathbb{R}}

\newcommand{\A}{\mathcal{A}}
\newcommand{\Aa}{\mathbb{A}}
\newcommand{\Ba}{\mathbb{B}}
\newcommand{\Ia}{\mathbb{I}}
\newcommand{\Ga}{\mathbb{M}_g}
\newcommand{\Fa}{\mathbb{F}}

\newcommand{\W}{\mathcal{W}}
\newcommand{\X}{\mathcal{X}}
\newcommand{\U}{\mathcal{U}}
\newcommand{\Z}{\mathcal{Z}}
\newcommand{\C}{\mathcal{C}}
\newcommand{\PC}{\mathcal{PC}}
\newcommand\norm[1]{\left\| #1\right\|}

\newcommand{\ds}{\displaystyle}

\newtheorem{definition}{\sc Definition}[section]
\newtheorem{teo}{\sc Theorem}[section]
\newtheorem{prop}{\sc Proposition}[section]

\newtheorem{lemma}{\sc Lemma}[section]

\newcommand{\dem}{\noindent {\bf Proof.} \mbox{}}
\newcommand{\findem}{\mbox{\framebox [1ex]{}}}
\begin{document}
\title[CONTROLLABILITY OF THE STRONGLY DAMPED WAVE EQUATION] { Controllability of the Strongly Damped Impulsive Semilinear Wave Equation with Memory and Delay}
%

\thanks{$\dagger$ This work has been supported by  Yachaytech University School of Mathematical Science and information technology, and Louisiana State University department of mathematics}
\author[C. Guevara and H. Leiva]{Cristi Guevara$^1$ and Hugo Leiva$^2$   }
\address{$^{1}$ Louisiana State University \\
          College of Science, Department of Mathematics \\
          Baton Rouge, LA 70803-USA} \email{cguevara@lsu.edu, cristi.guevara@asu.edu}
\address{$^{2}$ Yachaytech University \\
          School of Mathematical Science and information technology \\
          Hacienda San Jos\'e, San Miguel de Urcuqu\'i, Ecuador} \email{hleiva@yachaytech.edu.ec, hleiva@ula.ve}
\begin{abstract}
This article is devoted to study the interior approximated controllability of the strongly damped semilinear wave equation with memory, impulses and delay terms. The problem is challenging since the state equation contains memory  and impulsive terms yielding to potential unbounded control sequences steering the system to a neighborhood of the final state, thus fixed point theorems cannot be used directly. As alternative, the   A.E Bashirov and et al. techniques are applied and together with the delay allow the control solution to be directed to fixed curve in a short time interval and achieve our result.
\end{abstract}
\
\maketitle 
\section{Introduction}

The study of control systems has been of considerable interest for researchers motivated not only by engineering practices but also by biological process. Aiming to  improve manufacturing processes,  efficiency of energy use,  biomedical experimentation, diagnosis, robotics, biological control, systems among others. And surprisingly, it has become of great interest in the social, political and  economic spheres for understanding of the dynamics of business, social, and political systems.

Generally speaking,  control theory tackles how the  behavior of a systems can be modified by some feedback, specifically, how an arbitrary initial state can be directed, either exactly  or approximated close, to a given final state using a set of admissible controls. In addition, in practical control systems, abrupt changes, delays and dependance on prior behavior  are inherent  phenomena and they would  modify the controllability of the system.  Thus,  the conjecture is that  controllability of a system won't change due to  perturbations such as  delays, impulses  or some type of memories.

In this paper, we are concern with the approximated controllability  of strongly damped semilinear wave equation \eqref{eq:wave1} with memory, impulses and delay terms in $\Gamma =(0, \tau) \times  \Omega$, with  $\Omega \subseteq \R^{N} \,(N\geq1)$ a bounded domain and $\tau\in \R^+$,

\begin{equation}\label{eq:wave1}
\begin{split}
 w_{tt}+ \eta(-\Delta)^{1/2}w_{t}  + \gamma(- \Delta) w  = 1_{\varpi}u(t,x)+  \int_0^t M(t,s)g( w(s-r, x))ds\\
  + f(t,w(t-r),w_{t}(t-r),u(t,x)),
\end{split}\;
\end{equation}
 along with the initial-boundary conditions and impulses

\begin{equation}\label{initial}
\left\{\begin{array}{lcl}
 w(t,x)=0,&&\mbox{in}\; \partial  \Gamma=(0, \tau) \times \partial \Omega,\\
 \begin{split}
 &w(s,x)=\phi_1(s,x), \\
 &w_{t}(s,x)=\phi_2(s,x),
 \end{split} && \mbox{in} \; \Gamma_{r} = [-r,0] \times \Omega,\\
 w_{t}(t_{k}^{+},x) = w_{t}(t_{k}^{-},x)+I_{k}(t_k,w(t_{k},x),w_{t}(t_{k},x),u(t_{k},x)),& &  k=1,2,\dots,p,
  \end{array}
 \right.
\end{equation}
in the space $\Z^{1/2} = D\left( (-\Delta)^{1/2}\right) \times L^{2}(\Omega) $ as well as the distributed control $u \in L^{2} \left( [0,\tau], L^{2}(\Omega) \right)$ and $\Phi=(\phi_1,\phi_2) \in \mathcal{C}\left([-r,0]; \Z^{1/2} \right)$. Here, $\eta$, $\gamma$, $r$ are positive numbers, $r$ represents delay, $1_{\varpi}$ denotes the characteristic function on $\varpi$ an open nonempty subset of $\Omega$.

Assuming that the memory $M\in L^\infty \left( (0, \tau) \times  \Omega \right)$, the impulses $I_{k} :[0, \tau] \times \R^3 \rightarrow \R$ for  $k=1,2, \dots, p$,  the nonlinear functions  $g:\R \rightarrow \R$,   and $f:[0, \tau] \times \R^3 \rightarrow \R$ are smooth enough so that the problem \eqref{eq:wave1} admits only one mild solution on $[-r, \tau]$. Moreover,  for $u , w,v \in \R$ the following estimate holds
\begin{equation}\label{eq:hyp}
|f(t,w,v, u)| \leq a_{0}\sqrt{|w|^{2}+|v|^{2}} + b_{0}.
\end{equation}
Setting $J=[0,\tau]$ and  $J^{\prime} = [0,\tau] \backslash \{t_1, t_2, \dots, t_p \}$,  and let
\begin{align*}
\PC(J; \Z^{1/2} )  = &\left \{z: J  \rightarrow \Z^{1/2} \;/  z \in \C(J^{\prime}; \Z^{1/2} ),\; \exists\;  z(t_{k}^{+},\cdot),\; z(t_{k}^{-},\cdot)
 \;  \mbox{and} \; z(t_{k},\cdot)= z(t_{k}^{+},\cdot) \right\},
\end{align*}
equipped with the norm
$
\norm{z(\cdot)}_{_0} = \sup_{t \, \in  \,  [0,\tau]} \norm{z(t,\cdot )}_{ \Z^{1/2} }.
$ And for all  $z=(w,v)^{T} \in \Z^{1/2}$
$$
\norm{z}_{ \Z^{1/2} } ^2= \ds \int_{\Omega} \left(\norm{(-\Delta)^{1/2}w}^{2}+\norm{v }^{2}\right)dx.
$$

Equations of the form \eqref{eq:wave1} appear in a number of different contexts. One of them is in the study of motion of viscoelastic materials. For instance, in one-dimensional case, they model longitudinal vibration of a uniform, homogeneous bar with non-linear stress law. In two-dimension and three-dimension cases they describe antiplane shear motions of viscoelastic solids (see for details \cite{ FLO76, Knowles76, MM2000}). Furthermore, the strongly damped wave equation with memory has been used to  model  the deviation from the equilibrium configuration of a (homogeneous and isotropic) linearly viscoelastic solid  \cite{DL76}.   Note that including the impulses  on the system  \eqref{eq:wave1}, correspond to a pure mathematical interest and can physically interpret as  unexpected changes of state whose duration is negligible in comparison with the duration of the process.

In general, the existence of solutions for impulsive evolution equations with delays has been studied  by N. Abada, M. Benchohra and H. Hammouche \cite{ABH} and R.S. Jain and M.B. Dhakne in \cite{SB}. Moreover, the controllability of impulsive evolution equations is pretty well understood, just to mention some of the works, D.N. Chalishajar \cite{Chalishajar} studied the exact controllability of impulsive partial neutral functional differential equations with infinite delay and S. Selvi, M.M. Arjunan \cite{Selvi} studied the exact controllability for impulsive differential systems with finite delay.  L. Chen and G. Li \cite{Chen} studied the approximate controllability of impulsive differential equations with nonlocal conditions, using  measure of noncompactness and Monch fixed point theorem, and assuming that the nonlinear term $f(t,z)$ does not depend on the control variable. Adittionally, in a series of papers the approximate controllability of semilinear evolution equations with impulses is using the Rothe's Fixed Point Theorem  \cite{R3,R4,R5}, and using A.E. Bashirov and et at. technique which avoid fixed point theorems \cite{R6}

Motivated by the work of H. Larez, H.Leiva, J. Rebaza and A. Rios \cite{LLR, LLRR} on the approximate controlabiltiy of the semilinear strongly damped wave equation with and without impulses and our recent work on for the impulsive semilinear heat equation with memory and delay \cite{GL1}, we prove the interior approximate controllability of the strongly damped semilinear wave \eqref{eq:wave1} with memory, impulses and delay terms achieved by applying   A.E. Bashirov, N. Ghahramanlou, N. Mahmudov, N. Semi and H. Etikan technique \cite{B1,B2,B3}  avoiding fixed point theorems.

The structure of this paper is as follow: In section \ref{sec:formulation}, we present the abstract formulation of the strongly damped equation \eqref{eq:wave1}. Section \ref{lineal}, deals with the controllability of the linear problem. In section \ref{problem}, the approximated controllability of the strongly damped equation with memory, delay and impulses is proved.
\begin{definition} \label{def2}
{\rm (}{\sf Approximate Controllability}{\rm )} The system \eqref{eq:wave1} is said
to be approximately controllable on $J$, 
if for every $\Phi \in \C([-r,0]; \Z^{1/2} )$,
$z_1\in  \Z^{1/2} $, and $ \epsilon >0,$ there exists $u\in L^2\left(0,\tau; L^2(\Omega)\right)$ such
that the solution $z(t)$ of  \eqref{eq:wave1} corresponding to $u$
verifies:
 $$
 z(0) = \Phi(0) \ \ \mbox{and} \ \ \norm{z(\tau)-z_{1}}_{\Z^{1/2}}< \epsilon.
 $$
\end{definition}
The interior approximate controllability of the linear strongly damped wave equation
$$ 
w_{tt}+ \eta(-\Delta)^{1/2}w_{t} + \gamma(- \Delta) w  = 1_{\varpi}u(t,x),
$$
in $\Gamma$,   for all $\tau >0$, with initial-boundary  conditions
\begin{equation*}
\begin{array}{lll}
 w(t,x)=0, & & \mbox{in} \quad   \partial \Gamma,\\
 w(0,x)=w_0(x), \quad w_{t}(0,x)=w_1(x),& & \mbox{in} \quad  \Omega,
 \end{array}
\end{equation*}
was proven in  \cite{LLR}.

\section{Formulation of the Problem}\label{sec:formulation}
Let $\X = L^{2}(\Omega)= L^{2}(\Omega, \R)$, consider the linear
unbounded operator
\begin{equation*}
\begin{split}
\A: D(\A) \subset \X& \longrightarrow \X,\\
\phi& \longmapsto \A \phi= -\Delta \phi,
\end{split}
\end{equation*}
where $D(\A) = H^{2}(\Omega, \R) \cap H^{1}_{0}(\Omega, \R).$ Note that for  $\alpha \geq 0,$ the fractional powered spaces $\X^{\alpha}$ are given by
$$
\X^{\alpha} =D(\A^{\alpha}) = \left\{x \in \X : \sum_{n = 1}^{\infty} \lambda_{n}^{2 \alpha} \norm{ E_{n}x }^2 < \infty \right\}
$$
endowed with the norm
$$
\norm{x}_{\alpha}^2 = \|\A^{\alpha}x \|^2= \sum\limits_{n = 1}^{\infty}
\lambda_{n}^{2 \alpha} \| E_{n}x \|^2,
$$
where $\{ E_j \}$ is a family of complete orthogonal projections in $\X$; and  for  the Hilbert space $\Z^{\alpha} = \X^{\alpha} \times \X$ the corresponding norm  is
$$
 \norm{ \left( \begin{array}{c}w \\ v \end{array} \right)
 }_{ \Z^{\alpha}}^2 =  \|w \|^{2}_{r} + \|v \|^{2}.
$$
\begin{prop}\label{T1}
Given $j \geq 1$, the  operator $ P_j :  \Z^{\alpha} \rightarrow  \Z^{\alpha}$  defined by
\begin{equation}\label{eq:projection}
  P_j=\left[
\begin{array}{cc}
  E_j & 0 \\
  0 & E_j
\end{array}%
\right]
\end{equation}
is a continuous(bounded) orthogonal projections in the Hilbert space $\Z^{\alpha}$.
\end{prop}
Details and proofs can be found  in  \cite{LLR}.\\

Now, the system (\ref{eq:wave1})-(\ref{initial}) can be written as an abstract second order ordinary differential equation  in $\X$
\begin{equation}\label{eq: abstract0}
\left\{
\begin{array}{ll}
\begin{split}
w''+ \eta \A^{1/2}w' + \gamma \A w  =& b_{\varpi}u +\ds\int_0^tM(t,s)g^{e}(w(s-r))dr \\
&+ f^{e}(t,w(t-r),w'(t-r),u)
\end{split}, & t \in  (0, \tau]\; \mbox{and} \; t \neq t_{k},\\
 w(s,\cdot)=\phi_1(s,\cdot), \qquad w'(s,\cdot)=\phi_2 (s,\cdot),
 & s \in [-r,0], \\
w'(t_{k}^{+}) = w'(t_{k}^{-})+I^{e}_{k}(t_k, w(t_{k}),w'(t_{k}) ,u(t_{k})), & k=1,2, \dots, p,
\end{array}
\right.
\end{equation}
 where $\U=\X=L^{2}(\Omega)$,
\begin{eqnarray*}
I_{k}^{e}:&[0, \tau]\times \Z^{1/2} \times \U &\longrightarrow \qquad \X \\
&(t,w,v,u)(\cdot)&\longmapsto \quad I_{k}(t,w(\cdot),v(\cdot),u(\cdot)),
\end{eqnarray*}
\begin{eqnarray*}
f^{e}:&[0, \tau]\times \C(-r,0;  \Z^{1/2} ) \times \U & \longrightarrow \qquad \X\\
&(t,\phi_1,\phi_2,u)(\cdot)&\longmapsto \quad f(t,\phi_1(-r, \cdot),\phi_1(-r, \cdot),u(\cdot)),
\end{eqnarray*}
\begin{eqnarray*}
b_{\varpi}:&\U &\longrightarrow \qquad \U\\
&u(\cdot)&\longmapsto \quad 1_{\varpi}u(\cdot),
\end{eqnarray*}
and
\begin{eqnarray*}
g^{e}:&\C(-r,0;  \Z^{1/2} )  &\longrightarrow  \Z^{1/2} \\
&\Phi=\left(\begin{array}{c}
             \phi_1\\
             \phi_2
        \end{array}\right)&\longmapsto  g(\phi_1(\cdot-r)).
\end{eqnarray*}
 By changing variables $v=w'$, the second order
equation (\ref{eq: abstract0}) is written as a first order system of
ordinary differential equations with impulses and delay in the space $\Z^{1/2} =
\X^{1/2} \times \X$ as follows:
\begin{equation}\label{eq: abstract}
\left\{
\begin{array}{ll}
z' = \Aa z + \Ba_{\varpi}u +\ds\int_0^t\Ga(t,s,z(s-r))ds+ \Fa(t,z(t-r),u(t)),& z \in  \Z^{1/2}, \\
z(s)=\Phi(s), &s \in [-r,0],\\
z(t_{k}^{+}) = z(t_{k}^{-})+\Ia_{k}(t_k, z(t_{k}),u(t_{k})), & k=1,2, \dots, p,
\end{array}
\right.
\end{equation}
where $\Aa = \left(
   \begin{array}{rr}
     0 & I_{\X} \\ -\gamma \A & -\eta \A^{1/2}
   \end{array}\right)$ is a unbounded linear operator  with domain $D({\mathcal \Aa})= D(\A)
\times D(\A^{1/2})$, $I_{\X}$ is the identity in $\X, $
$ z = \left(
  \begin{array}{c}
    w \\ v
  \end{array}
  \right),$
$u\in \C\left([0,\tau];\U\right)$,
   $\Phi=\left(\begin{array}{c}
             \phi_1\\
             \phi_2
        \end{array}\right) \in \C\left(-r,0;  \Z^{1/2} \right),$
\begin{eqnarray*}
\Ba_{\varpi}:& \U& \longrightarrow \qquad \Z^{1/2} \\
&u&\longmapsto \quad
\left(\begin{array}{c}
             0\\
             b_{\varpi}u
        \end{array}\right),
\end{eqnarray*}
\begin{eqnarray*}
\Ia_{k}:&[0, \tau]\times \Z^{1/2} \times \U& \longrightarrow \qquad \Z^{1/2} \\
&(t, z,u)&\longmapsto \quad
\left(\begin{array}{c}
             0\\
             I_{k}^{e}(t,w,v,u)
        \end{array}\right),
\end{eqnarray*}

\begin{eqnarray}\label{functionF}
\Fa:& [0, \tau] \times \C(-r,0;  \Z^{1/2} ) \times \U & \longrightarrow \qquad \Z^{1/2}\\
&(t, \Phi,u)&\longmapsto \quad
\left(
  \begin{array}{c}
    0 \\f^{e}(t,\phi_1(-r),\phi_2(-r),u)
  \end{array}
  \right),\notag
  \end{eqnarray}
and
\begin{eqnarray*}
\Ga:&[0,\tau]\times [0,\tau]\times \C(-r,0;  \Z^{1/2} )  &\longrightarrow  \Z^{1/2} \\
&(t,s,\Phi)&\longmapsto \left(
  \begin{array}{c}
    0 \\M(t,s) g^{e}(\Phi)
  \end{array}
  \right).
\end{eqnarray*}
The continuous inclusion $\X^{1/2} \subset \X,$ together with  the condition \eqref{eq:hyp}, yields the  following proposition
\begin{prop} The function $\Fa$  as defined by \eqref{functionF} satisfies
\begin{equation}\label{eq: boundF}
\norm{\Fa(t, \Phi, u)}_{ \Z^{1/2}} \leq \tilde{a}\norm{\Phi(\cdot-r)}_{ \Z^{1/2} }+\tilde{b},
\end{equation}
 for all $(t, \Phi, u) \in [0, \tau]\times \C(-r,0;  \Z^{1/2} ) \times \U$ and $\tilde{a},\tilde{b}>0.$
\end{prop}
It has been  proved in  \cite{chen1, Chen}  that the operator $\Aa$ generates a strongly continuous and  analytic semigroup
$\left\{ T(t) \right\}_{t \geq 0}$ in the space $\Z^{1/2} = \X^{1/2} \times \X$. Furthermore,  Lemma 2.1 in \cite{HL} yields
\begin{prop}\label{T2}
The semigroup $\left\{ T(t) \right\} _{t \geq 0}$, generated by the operator $\Aa$ is compact and  represented by
\begin{equation}\label{damp}
  T(t)z =\sum_{j=1}^{\infty}e^{A_{j}t}P_jz,
  \end{equation}
  with $z\in \Z^{1/2},\,$ $t \geq 0$, here $\left\{ P_j\right\} _{j \geq 1}$ is a complete family of
orthogonal projections \eqref{eq:projection} in the Hilbert space $ \Z^{1/2} $ and $A_j = R_{j}P_{j}$ for $j=1,2, \cdots,$ where
\begin{equation*}
R_j=\left(
  \begin{array}{cc}
     0 &  1 \\
     -\gamma \lambda _j &  -\eta \lambda^{1/2}_{j}
  \end{array}
  \right),
  \end{equation*}
  with eigenvalues $\lambda = -\lambda^{1/2}_{j} \left( \frac{\eta \underline{+}
\sqrt{\eta^{2} - 4\gamma }} {2} \right).$  Moreover, $e^{A_jt}=e^{R_jt}P_j$,
and $ A_j^*=R_j^*P_j$
$$
R^*_j=\left(\begin{array}{cc}
             0 & -1\\
            \gamma \lambda _j &  -\eta \lambda^{1/2}_{j}
            \end{array}
            \right).
$$
\end{prop}

\section{Controllability of the Linear System} \label{lineal}
In this section we present some characterization of the interior approximate controllability of the linear
strongly damped  wave equations.

The initial value problem
\begin{equation}\label{eq:Linear}
\left\{\begin{array}{ll}
             z^{\prime}= \Aa z + \Ba_{\varpi}u, \\
             z(t_0)=z_0,
        \end{array}\right.
\end{equation}
with $z, z_0\in  \Z^{1/2} $ and $u\in L^{2}(0,\tau;\U)$,
admits only one mild solution for $ t \in [t_0,\tau]$, given by
\begin{equation}\label{eq: mildlinear}
z(t)=T(t-t_0)z_0+\int_{t_0}^{t}T(t-s)\Ba_{\varpi}u(s)ds.
\end{equation}

\begin{definition}\label{def3}
For system \eqref{eq:Linear} and $\tau>0,$
the controllability map is given by
\begin{equation} \label{eq:control1}
\begin{split}
G_{\tau\delta}: L^2(\tau-\delta,\tau;\U) \longrightarrow& \Z^{1/2}\\
u\longmapsto& \int_{\tau-\delta}^{\tau}T(\tau-s)\Ba_{\varpi}u(s)ds,
\end{split}
\end{equation}
 and its adjoint by
\begin{equation*}
 \begin{split}
G^*_{\tau\delta}:  \Z^{1/2} \longrightarrow& L^2(\tau-\delta,\tau;\U)\\
z\longmapsto& \Ba_{\varpi}^{*}T^{*}(\tau-\cdot)z.
\end{split}
\end{equation*}
 The Gramian controllability operator is
\begin{equation}\label{Q2}
Q_{\tau \delta} = G_{\tau\delta}G_{\tau\delta}^{*}= \int_{\tau-\delta}^{\tau}T(\tau-t)\Ba_{\varpi}\Ba_{\varpi}^{*}T^{*}(\tau-t)dt.
\end{equation}
\end{definition}

\begin{lemma}\label{l:propertiesG}
The following statements are equivalent to the approximate controllability of the linear system \eqref{eq:Linear}  on $[\tau-\delta,\tau]$.

\begin{enumerate}
\item $\overline{\textup{Rang}(G_{\tau\delta})}=\Z^{1/2}$.

\item $Ker(G_{\tau\delta}^{*})=\{0\}$.
\item If $0 \neq z\in \Z^{1/2}$, then $\langle Q_{\tau\delta}z,z\rangle>0$.

\item  $\ds\lim_{\alpha \to 0^+} \alpha(\alpha I +Q_{\tau\delta})^{-1}z =0 $.
\item For all $z \in \Z$, $0<\alpha\leq1$  and $u_{\alpha}= G_{\tau\delta}^{*}(\alpha I+Q_{\tau\delta})^{-1}z$ we have that
$$
G_{\tau\delta}u_{\alpha}= z -\alpha(\alpha I+Q_{\tau\delta})^{-1}z
$$
Hence, 
$$
\lim_{\alpha \to 0} G_{\tau\delta}u_{\alpha}= z
$$ 
and error $E_{\tau\delta}z$ of this approximation is
$$
E_{\tau\delta}z=\alpha(\alpha I+ Q_{\tau\delta})^{-1}z.
$$
In addition, for each $v \in L^2(\tau-\delta,\tau;\U)$,  the sequence of controls
$$
u_{\alpha}= G_{\tau\delta}^{*}(\alpha I+Q_{\tau\delta})^{-1}z+(v-G_{\tau\delta}^{*}(\alpha I+Q_{\tau\delta})^{-1}G_{\tau\delta}v),
$$
satisfies
$$
G_{\tau\delta}u_{\alpha}= z -\alpha(\alpha I+Q_{T\delta})^{-1}(z+G_{\tau\delta}v)
$$
and
$$
\lim_{\alpha \to 0} G_{\tau\delta}u_{\alpha}= z,
$$
with  error $E_{\tau\delta}z=\alpha(\alpha I+ Q_{\tau\delta})^{-1}(z+G_{\tau\delta}v)$.

\end{enumerate}
\end{lemma}
In general, for a pair $\W$ and $\Z$ of Hilbert spaces,  a linear bounded operator $G:\W \rightarrow \Z$ with dense range  Lemma \ref{l:propertiesG} also holds and the proof  can be found in  \cite{E, F, CP1, CP2, LMS2}.\\
This Lemma implies that for $0 < \alpha \leq 1$, the family of linear operators $\Gamma_{\alpha\tau\delta} :\Z^{1/2} \rightarrow L^2(\tau-\delta,\tau;\U)$, defined by
\begin{equation}\label{rappro}
\Gamma_{\alpha\tau\delta}z=  G_{\tau\delta}^{*}(\alpha I+ Q_{\tau\delta})^{-1}z,
\end{equation}
is a approximated right inverse operator of $G_{\tau\delta}$, in the sense that
\begin{equation}
\lim_{\alpha \to 0} G_{\tau\delta} \Gamma_{\alpha\tau\delta}= I.
\end{equation}
in the strong topology.

\begin{lemma}\label{lema}
$Q_{\tau\delta} > 0$, if and only if, the linear system \eqref{eq:Linear} is controllable on $[\tau-\delta, \tau]$.
Moreover, given an initial state $y_0$ and a final state $z^{1}$ we can find a sequence of controls $\{u_{\alpha}^{\delta}\}_{0 <\alpha \leq 1} \subset L^2(\tau-\delta,\tau;\U)$
$$
u_{\alpha}=u_{\alpha ,\delta}= G_{\tau\delta}^{*}(\alpha I+ G_{\tau\delta}G_{\tau\delta}^{*})^{-1}(z^{1} - T(\tau)y_0),
$$
such that the solutions $y(t)=y(t,\tau-\delta, y_0, u_{\alpha})$ of the initial value problem
\begin{equation}\label{IVL}
\left\{
\begin{array}{l}
y'=\Aa y+\Ba u_{\alpha}(t), \ \  y \in Z, \ \ t>0,\\
y(\tau-\delta) = y_0,
\end{array}
\right.
\end{equation}
satisfies
$$
\lim_{\alpha \to 0^{+}}y(\tau,\tau-\delta, y_0, u_{\alpha}) = z^{1}.
$$
that is,
$$
\lim_{\alpha \to 0^{+}}y(\tau)= \lim_{\alpha \to 0^{+}}\left\{T(\delta)y_0 + \int_{\tau-\delta}^{\tau}T(\tau-s)\Ba u_{\alpha}(s)ds \right\}= z^{1}.
$$
\end{lemma}

\section{Controllability of the semilinear System}\label{problem}
In this section we prove the interior approximate controllability of
the semilinear strongly damped wave equation with memory, impulses and delay, 
Theorem \ref{main}. The idea behind of the proof is that for a given final state $z_1,$  a sequence of controls is constructed, so that the initial condition $\Phi$ is steering to a small ball around $z_1.$ All of this is achieved thanks to the delay, which allows to pullback the family of controls solutions to a fixed trajectory in short time interval, as illustrated below:
\begin{center}
\includegraphics[height=2.5 in]{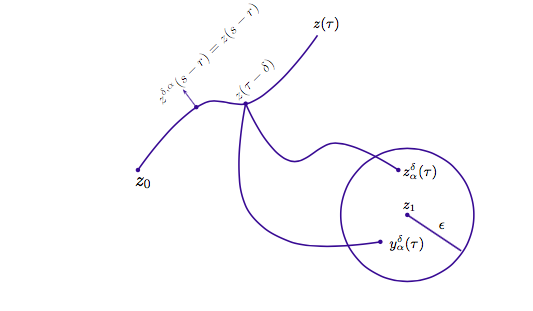}
\end{center}
\begin{teo}\label{main}  The wave equation  \eqref{eq:wave1} under \eqref{initial} and \eqref{eq:hyp} is approximately
controllable on $[0,\tau]$ with $\tau>0$.
\end{teo}

\dem Provided $\epsilon>0$, $\Phi\in \mathcal{C}$ and a final state $z_{1}$, we are aiming to find a control $u_{\alpha}^{\delta}\in L^{2}([0,\tau];\U)$ such that it evolves the system \eqref{eq:wave1}  from  $\Phi(0)$ to an $\epsilon-$ ball around  $z_{1}$ on time $\tau$. That is,
considering the  the abstract formulation for  \eqref{eq:wave1}  discussed in section,
 recall that for all $\Phi \in \C(-r,0;  \Z^{1/2} )$ and $u\in \C (0,\tau;\U)$ the   system \eqref{eq: abstract} given by
\begin{equation*}\label{formulation}
\left\{
\begin{array}{ll}
z' = \Aa z + \Ba_{\varpi}u +\ds\int_0^t \Ga(t,s,z(s-r))ds+ \Fa(t,z(t-r),u(t)),& z \in  \Z^{1/2}, \\
z(s)=\Phi(s), &s \in [-r,0],\\
z(t_{k}^{+}) = z(t_{k}^{-})+\Ia_{k}(t_k, z(t_{k}),u(t_{k})), & k=1,2,\dots,p.
\end{array}
\right.
\end{equation*}
We want to prove that, for $\alpha>0$ and $0<\delta<\min\{\tau-t_{p},r\}$, there exists control $u_{\alpha}^{\delta}\in L^{2}([0,\tau];\U)$ such that corresponding of solutions $z^{\delta \alpha}$ of  \eqref{eq: abstract} satisfies:
$$
\parallel z^{\delta \alpha}(\tau)-z_{1}\parallel\leq \epsilon,
$$
and
\begin{equation}\label{eq: mild}
\begin{split}
z^{\delta \alpha}(t) =&  T(t)\Phi(0)+\int_{0}^{t}T(t-s)\left[\Ba_{\varpi}u(s)+ \int_0^s \Ga(t,s,z(l-r))dl\right]ds+ \\
&+  \int_{0}^{t}T(t-s)\Fa(s,z(s-r),u(s))ds+
\sum_{0 < t_k < t} T(t-t_k )\Ia_{k}(t_k,z(t_k), u(t_k)),
\end{split}
\end{equation}
for $0\leq t \leq\tau,$ admitted by  the system \eqref{eq: abstract}.

Consequently,  take $u\in L^{2}([0,\tau];\U)$ and its corresponding solution $z(t)=z(t,0,\Phi,u)$ of the initial value problem \eqref{eq: abstract} and define the controls $u_{\alpha}^{\delta}\in L^{2}([0,\tau];\U)$ for $\alpha\in (0,1]$, such that
$$
u_{\alpha}^{\delta}(t)=\left\{\begin{array}{lc}
                         u(t),& 0\leq t\leq \tau-\delta, \\
                         u_{\alpha}(t)= \Ba^{*}T^{*}(\tau-t)(\alpha I+ G_{\tau\delta}G_{\tau\delta}^{*})^{-1}(z_{1} - T(\delta)z(\tau-\delta)), &\tau-\delta\leq t\leq \tau,
                       \end{array}\right.
$$
Since $0<\delta<\tau-t_{p}$,  then the corresponding solution $z^{\delta,\alpha}(t)=z(t,0,\Phi,u_{\alpha}^{\delta})$ of the initial value problem at time $\tau$ is written by:
\begin{align*}
z^{\delta,\alpha}(\tau) =&
 T(\tau)\Phi(0)
+\int_{0}^{\tau}T(\tau-s)
	\left[ \Ba_{\varpi} u_{\alpha}^{\delta} (s)
		+ \int_0^s \Ga(z^{\delta,\alpha}(l-r))dl\right]ds+ \\
	&+  \int_{0}^{\tau}T(\tau-s)\Fa(s,z^{\delta,\alpha}(s-r),u_{\alpha}^{\delta}(s))ds+
\sum_{0 < t_k < \tau} T(t-t_k )\Ia_{k}(t_k,z^{\delta,\alpha}(t_k), u_{\alpha}^{\delta}(t_k))\\
=&T(\delta)\left\{T(\tau-\delta)\Phi(0)
+\int_{0}^{\tau-\delta}T(\tau-\delta-s) \left(\Ba_{\varpi} u_{\alpha}^{\delta} (s)+\Fa(s,z^{\delta,\alpha}(s-r),u_{\alpha}^{\delta}(s))\right)ds\right.\\
&\qquad+\int_{0}^{\tau-\delta}T(\tau-\delta-s)  \int_0^s \Ga(s,l, z^{\delta,\alpha}(l-r))dlds\\
&\qquad\left.+ \sum_{0 < t_k < \tau-\delta} T(t-\delta-t_k )\Ia_{k}(t_k,z^{\delta,\alpha}(t_k), u_{\alpha}^{\delta}(t_k))\right\}+\\
& + \int_{\tau-\delta}^{\tau}T(\tau-s)\left(\Ba u_{\alpha}(s)+
\Fa(s,z^{\delta,\alpha}(s-r),u_{\alpha}^{\delta}(s))+\int_0^s\Ga(s,l,z^{\delta,\alpha}(l-r))dl\right)ds,
\end{align*}
hence,
\begin{align*}
z^{\delta,\alpha}(\tau)  = &T(\delta)z(\tau-\delta)+ \int_{\tau-\delta}^{\tau}T(\tau-s)\left(\Ba u_{\alpha}(s)+
\Fa(s,z^{\delta,\alpha}(s-r),u_{\alpha}^{\delta}(s))\right)ds \\
&+ \int_{\tau-\delta}^{\tau}T(\tau-s)\int_0^s\Ga(s,l,z^{\delta,\alpha}(l-r))dlds.
\end{align*}
Since the corresponding solution $y^{\delta,\alpha}(t)=y(t,\tau-\delta,z(\tau-\delta),u_{\alpha})$ of the linear initial value problem \eqref{eq: mildlinear} at time $\tau$ is
$$
y^{\delta,\alpha}(\tau)=T(\delta)z(\tau-\delta)+ \int_{\tau-\delta}^{\tau}T(\tau-s)\Ba u_{\alpha}(s)ds,
$$
then
$$
z^{\delta,\alpha}(\tau)-y^{\delta,\alpha}(\tau)=
\int_{\tau-\delta}^{\tau}T(\tau-s)\left(\int_{0}^{s}\Fa(s,z^{\delta,\alpha}(s-r),u_{\alpha}^{\delta}(s)+\Ga(s,l,z^{\delta,\alpha}(l-r))dl)\right)ds.
$$
As consequence, the following estimate yields
\begin{align*}
  \norm{ z^{\delta,\alpha}(\tau)-y^{\delta,\alpha}(\tau)} & \leq \int_{\tau-\delta}^{\tau} \norm{ T(\tau-s)}\left( \tilde{a}\norm{\Phi(s-r)}+\tilde{b}\right)ds \\
   & +  \int_{\tau-\delta}^{\tau}\norm{ T(\tau-s)}\int_{0}^{s}\norm{\Ga(s,l,z^{\delta,\alpha}(l-r))}dlds.
\end{align*}
Notice that $0< \delta< r$ and $\tau-\delta \leq s\leq \tau$, thus $$l-r \leq s-r \leq \tau-r< \tau-\delta,$$
hence
$$
z^{\delta,\alpha}(l-r)=z(l-r) \quad\mbox{and}\quad z^{\delta,\alpha}(s-r)=z(s-r),
$$
which implies that there exists $\delta>0$ such that
\begin{align*}
 \norm{z^{\delta,\alpha}(\tau)-y^{\delta,\alpha}(\tau)} & \leq \int_{\tau-\delta}^{\tau}\norm{ T(\tau-s)}\left( \tilde{a}\norm{z(s-r)}+\tilde{b}\right)ds \\
  &\quad +  \int_{\tau-\delta}^{\tau}\norm{T(\tau-s)}\int_{0}^{s}\norm{ \Ga(s,l,z(l-r))} dlds  \\
   & <  \displaystyle\frac{\epsilon}{2}.
\end{align*}
Additionally, by Lemma \ref{lema} we can chose $\alpha>0$ such that
$$
\begin{array}{lll}
 \norm{ z^{\delta,\alpha}(\tau)-z_{1}}  & \leq &  \norm{ z^{\delta,\alpha}(\tau)-y^{\delta,\alpha}(\tau)} +  \norm{ y^{\delta,\alpha}(\tau)-z_{1}}  <  \frac{\epsilon}{2}+ \frac{\epsilon}{2}=\epsilon,
\end{array}
$$
which completes our proof.

\hfill \findem

\section{Final Remarks}\label{final}
\noindent

The technique employed here can be applied for  the   impulsive semilinear beam equation with memory and delay in a bounded domain $\Omega \subseteq \R^{N} \,(N\geq1)$,:
\begin{equation}
\left\{%
\begin{array}{lll}
  \displaystyle{\partial^{2} y(t,x) \over \partial t^{2}}
  & = &  2\beta\Delta\displaystyle\frac{\partial y(t,x)}{\partial t}- \Delta^{2}y(t,x) + \int_{0}^{t}M(t-s)g(y(s-r,x))ds\\
   & &  + u(t,x) + f(t,y(t-r),y_{t}(t-r),u),\; \mbox{in}\; \Omega_{\tau},\;t\neq t_{k}, \\
  y(t,x) & = & \Delta y(t,x)= 0 , \ \ \mbox{on}\; \Sigma_{\tau}, \\
  y(s,x) & = & \phi(s,x), \ \ y_{t}(s,x)=\psi(s,x),\quad \mbox{in}\;\Omega_{-r},\\
y_{t}(t_{k}^{+},x) & = & y_{t}(t_{k}^{-},x) +I_{k}(t_{k},y(t_{k},x), y_{t}(t_{k},x),u(t_{k},x)), \quad k=1, \dots, p,
\end{array}%
\right.
\end{equation}
where  $\Omega_{\tau}=(0,\tau]\times \Omega$, $\Sigma_{\tau}=(0,\tau)\times \partial \Omega$, $\Omega_{-r}=[-r,0]\times \Omega$.

Furthermore, we believe that the same technique  can be applied for  controlling diffusion processes systems involving compact semigroups. In particular, our result can be formulated in a more general setting  for the semilinear evolution equation with impulses, delay and memory  in a Hilbert space $\Z$
\begin{equation*}\label{eq:formulation}
\left\{
\begin{array}{ll}
z' = \Aa z + \Ba_{\varpi}u +\ds\int_0^t \Ga(t,s,z(s-r))ds+ \Fa(t,z(t-r),u(t)),& z \in  \Z \\
z(s)=\Phi(s), &s \in [-r,0]\\
z(t_{k}^{+}) = z(t_{k}^{-})+\Ia_{k}(t_k, z(t_{k}),u(t_{k})), &
\end{array}
\right.
\end{equation*}
where $u\in \C ([0,\tau];\U)$, $\U$ is another Hilbert space, $B :\U \longrightarrow \Z$ is a bounded linear operator,$J_{k}, F:[0, \tau]\times C(-r,0; \Z) \times \U \rightarrow \Z$, $\Aa :D(\Aa) \subset \Z \rightarrow \Z$ is an unbounded linear operator in $\Z$ that generates a strongly continuous semigroup according to
Lemma 2.1 from \cite{HL}:
\begin{equation}\label{damp2}
  T(t)z =\sum_{nj=1}^{\infty}e^{A_{j}t}P_jz%
  \mbox{, } \ \ z\in \Z \mbox{, } \ \ t \geq 0,     
\end{equation}
where  $\left\{ P_j\right\} _{j \geq 0}$ is a complete family of orthogonal projections in the Hilbert space $\Z$ and
\begin{equation}
\|F(t,\Phi,u) \|_{\Z}  \leq   \tilde{a} \|\Phi(-r)\|_{\Z}^{\alpha_{0}} +\tilde{b}, \quad \forall (t, \Phi, u) \in [0, \tau]\times C(-r,0;  \Z^{1/2} ) \times \U.
\end{equation}

%

\end{document}